\documentclass [nomath,noamsfonts]{amsart}

\newcommand{\force}{\mbox{$\Vdash$}}

\newcommand{\ga}{\alpha}
\newcommand{\gb}{\beta}
\newcommand{\grg}{\gamma}
\newcommand{\gd}{\delta}
\newcommand{\gre}{\varepsilon}

\newcommand{\gz}{\zeta}

\newcommand{\gk}{\kappa}
\newcommand{\gl}{\lambda}
\newcommand{\gm}{\mu}
\newcommand{\gn}{\nu}

\newcommand{\gr}{\rho}
\newcommand{\gs}{\sigma}

\newcommand{\go}{\omega}

\newcommand{\gP}{\Pi}
\newcommand{\gS}{\Sigma}

\newcommand{\bgS}{\boldsymbol \Sigma}

\newcommand{\ha}{\aleph}

\theoremstyle{plain}
\newtheorem{theorem}{Theorem}[section]
\newtheorem{lemma}[theorem]{Lemma}
\newtheorem{prop}{Proposition}[section]

\theoremstyle{definition}
\newtheorem{definition}{Definition}[section]
\theoremstyle{remark}
\newtheorem{coro}{Corollary}[section]

\newcommand{\beq}{\begin{equation}}
\newcommand{\eeq}{\end{equation}}
\newcommand{\bds}{\begin {itemize}}
\newcommand{\eds}{\end {itemize}}
\newcommand{\bdf}{\begin{definition}}
\newcommand{\blm}{\begin{lemma}}
\newcommand{\edf}{\end{definition}}
\newcommand{\elm}{\end{lemma}}
\newcommand{\bthm}{\begin{theorem}}
\newcommand{\ethm}{\end{theorem}}
\newcommand{\bprp}{\begin{prop}}
\newcommand{\eprp}{\end{prop}}
\newcommand{\bcr}{\begin{coro}}
\newcommand{\ecr}{\end{coro}}

\newcommand{\rarrow}{{\rightarrow}}

\newcommand{\lrarrow}{{\leftrightarrow}}

\input {fullpage.sty}
\usepackage{amssymb}
\usepackage{amsthm}
\usepackage{amsbsy}
\newcommand{\coll}{\hbox{\ Coll \ }}
\newcommand{\mP}{\hbox{{\bf P}}}
\newcommand{\Hom}{\mbox{$(H_{\go_1},\gre)$}}
\newcommand{\bprf}{\hbox{{\bf Proof :}}}
\newcommand{\eprf}{\mbox{$\qed \\$}}
\newcommand{\mPeps}{\mbox{$\bf P_{\gre}$}}
\newcommand{\tlgs}{\mbox{$\tilde \gs$}}
\newcommand{\tlh}{\mbox{$\tilde h$}}
\numberwithin{equation}{section}
\begin {document}
\title [The definable tree property]{On the consistency of the definable
\\ tree property on $\aleph_1$}
\author {Amir Leshem}
\address{Institute of Mathematics, Hebrew University, Jerusalem
91904, Israel}
\address{Circuit and systems group, facaulty of information technology and systems ,
Delft University of Technology, Mekelweg 4, 2628CD Delft, The Netherlands}
\email{leshem@cas.et.tudelft.nl}
\email{leshem@math.huji.ac.il}
\date {\today}
\begin{abstract}
In this paper we prove the equiconsistency of
``Every $\go_1-$tree which is
first order definable over $\Hom$ has a cofinal branch'' with the existence
of a $\Pi^1_1$ reflecting cardinal. We also prove that the addition of
$MA$ to the definable tree property increases the consistency strength to 
that of
a weakly compact cardinal. Finally we comment on the generalization to
higher cardinals.
\end{abstract}
\maketitle

\bibliographystyle{plain}
\section{Introduction}
A well known result of Aronszajn is the existence of an Aronszajn tree
on $\ha_1$, i.e. there exist an $\go_1-$tree with no uncountable branch.
The construction uses the axiom of choice and therefore does not give a
definable such tree.

Sierpinski \cite{Si} and Kurepa \cite{Kr} proved
that if
Ramsey theorem holds for $\gk$ then $\gk$ is a
strong limit cardinal. Erd\H{o}s \cite{ErTa} proved that such a
$\gk$ is inaccessible. They also provided counterexamples to
Ramsey theorem on small cardinals. These
counterexamples explicitly used a well-ordering of
$P(\gk)$. The proof  raises
the question whether we can find a definable
counterexample.

In a straightforward generalization of Aronszajn's proof, Specker proved
the existence
of Aronszajn trees for every succesor $\gk^+$ s.t
$\gk^{<\gk}=\gk$. This raised the question whether the
GCH was needed for this result.
Mitchell and Silver \cite{Mi} have proved that the tree
property
on $\ha_2$ is equiconsistent with the existence of a weakly
compact cardinal. Magidor and Shelah \cite{MgSh} have proved the
consistency of the tree property on $\ha_{\go+1}$
 from the existence of very large cardinals.

Mitchell's forcing for the tree property works well if one
tries to obtain the tree property on two non-consecutive
cardinals (e.g. $\ha_2$ and $\ha_4$). However, his methods
fail to prove the consistency of the tree property on $\ha_2$
and $\ha_3$ together. By a result of Magidor this is indeed a
substential difficulty since the consistency strength of the
tree property on $\ha_2$ and $\ha_3$ together is much higher
then a weakly compact (e.g. it implies the existence of
$0^{\sharp}$). Abraham \cite{Ab} has proved the consistency of the tree
property on both $\ha_2$ and $\ha_3$ from a supercompact
cardinal and a weakly compact cardinal above it. This result
was extended by Cummings and Foremann \cite{CuFo} which proved
the consistency
of the tree property on all $\ha_n$'s from many supercompacts.

Kunen, and Shelah and Harrington \cite{HaSh} considered the consistency strength of
obtaining Lebesgue measurability of projective sets together with Martin's
axiom. They proved the equiconsistency of this theory with the existence
of weakly compact cardinals. This shows that adding Martin's axiom to
projective measurability increases the consistency strength.

In this paper we consider the consistency strength of the definable tree
property on $\ha_1$, i.e. the existence of a model in which every
$\go_1-$tree which is first order definable (with parameters)
over $\Hom$, has a cofinal branch. Our proof method answers
also a related question regarding definable
counterexamples to Ramsey theorem
on $\ha_1$.

In section (\ref{tree_def}) we define the exact meaning of a definable
$\gk-$tree, and study some variations used in this paper.
In section (\ref{reflect}) we define the $\Pi^1_1$ reflecting cardinals
and derive an extension property which is similar to the extension
property of weakly compact cardinal. We also bound the consistency
strength of the existence of a $\Pi^1_1$ reflecting cardinal by the
existence of a Mahlo cardinal.
In section (\ref{Forcing_model}) we show that by forcing with
the well known Levy collapse of a $\Pi_1^1$ reflecting cardinal $\gk$
to $\ha_1$, we obtain a
model of the definable tree property. We also prove that in this model
a definable Ramsey theorem on $\ha_1$ holds.
In section (\ref{lower_bound})
we show that our assumptions on $\gk$ are necessary by proving that
if the definable tree property on $\ha_1$ holds, then
\[
L \models \ha_1 \hbox{\ is a $\Pi^1_1$ reflecting cardinal}.
\]
In section (\ref{DTPMA}) we prove the equiconsistency of the definable
tree property on $\ha_1+$  MA with the existence of a weakly
compact cardinal. This is done by exploiting the methods of
Kunen and Shelah and Harrington \cite{HaSh}. We also comment that adding
any reasonable failure of GCH does not add to the consistency strength.
Finally in section (\ref{DTPGCH}) we comment on the consistency of the
definable tree property on higher cardinals with GCH. We describe a
forcing to get the definable tree property on all $\ha_n$'s together
with GCH by using $\go$ many $\Pi^1_1$ reflecting cardinals.
\section{Definable $\gk-$trees}
\label{tree_def}
In this section we define various notions of definable $\gk-$trees
and study the relationship between these notions.
First define the usual notion of a $\gk-$tree.
\bdf

\bds
\item A tree is a partially ordered set$(T,<_T)$ such that
for any $t \in T$ the set $\{s \in T |s<_Tt \}$ of predecessors of $t$
is well-ordered under $<_T$, and there is a root $r \in T$ such that for any
$t \in T$, such that $t \not = r$, $r <_T t$.
\item The $\ga$'th level of $T$ denoted by $T_{\ga}$ is the set of
elements of T whose set of $T-$predecessors has order type $\ga$.
\item A tree $(T,<_T)$ is a $\gk-$tree if $|T|=\gk$ and for every $\ga$
$|T_{\ga}|<\gk$.
\item A branch is a $<_T$ linearly ordered subset of $T$.
\item A cofinal branch is a branch which intersects every level of $T$.
\eds
\edf
Next we would like to give a definition of a definable $\gk-$tree.
Three different notions naturally arise
\bdf

\bds
\item A $\gk-$tree is \em{definable in the strict sense} if its
underlying set is $\gk$, and $<_T$ is
$\bgS_{\go}\left((H_{\gk},\in)\right)$.
\item A $\gk-$tree is \em{definable in the wide sense} if its
underlying set $T$ and $<_T$ are both
$\bgS_{\go}\left((H_{\gk},\in)\right)$ and $T$ has definable cardinality
$\gk$, i.e. there is a bijection
$f:\gk \lrarrow T$ which is  $\bgS_{\go}\left((H_{\gk},\in)\right)$.
\item  A $\gk-$tree is \em{definable in the very wide sense} if its
underlying set $T$ and $<_T$ are both
$\bgS_{\go}\left((H_{\gk},\in)\right)$.
\eds
\edf
 Obviously being definable in the wide sense
implies being definable in the very wide sense, and
being definable in the strict sense
implies being definable in the wide sense.
The following proposition states that being definable in the wide sense
is almost equivalent to being definable in the strict sense.
\bprp
If $(T,<_T)$ is a $\gk-$tree definable in the wide sense then there is a
tree $(\gk,<_{T'})$ isomorphic to $(T,<_T)$ which is definable in
the strict sense.
\eprp
\bprf
Let $f$ be the definable bijection $f:\gk \rarrow T$, and let
$\psi(x,y,z)$ be the definition of $<_T$ from the parameter $z$.
Now define $<_{T'}$ by $\ga<_{T'}\gb$ iff $\psi(f(\ga),f(\gb),z)$. The rest
trivially follows.
\eprf

Using this proposition we will not distinguish between the strict and the
wide notions of definability. Also notice that if there is a
well-ordering of $H_{\gk}$ which is
$\bgS_{\go}\left((H_{\gk},\in)\right)$ then the last two notions
coincide. However in general this is not the case, and we will
distinguish between definability in the very wide sense and definability
in the strict sense which we will adopt as the definition of definability.
So from now on, a definable $\gk-$tree is a
tree definable in the strict sense.
\section{$\Pi^1_1$ reflecting cardinals}
\label{reflect}

Let $\gk$ be a cardinal. We say that $\gk$ is $\Pi^m_n$ reflecting, if
$\gk$ is inaccessible and
for every $A \subseteq V_{\gk}$ {\em definable} over $V_{\gk}$ (with
parameters) and for
every $\Pi^m_n$ sentence
$\Phi$, such that
\[
(V_{\gk},\gre,A) \models \Phi
\]
there is an $\ga<\gk$ such that
\[
(V_{\ga},\gre,A \cap V_{\ga})  \models \Phi.
\]
The $\gP^m_n$ reflecting cardinals are a lightface analog of the $\gP^m_n$
indescribable cardinals. However they have a much weaker consistency strength.
If $\gk$ is $\Pi^1_1$ reflecting, then an easy consequence is the fact that
$\gk$
is the $\gk 'th$ inaccessible cardinal, since inaccessibility is a $\Pi^1_1$
property, and for each $\ga<\gk$ being the $\ga 'th$ inaccessible is also
a $\Pi^1_1$ property.
Next we prove the following lemma:
\blm
Suppose $\gk$ is a Mahlo cardinal then for every $m,n$,\\
$\{\ga<\gk | \ga \hbox{\ is a $\Pi^m_n$ reflecting cardinal} \}$
is stationary.
\elm
\bprf
Let $\gk$ be Mahlo, and let $C$ be a club. We shall find a $\gm \in C$
such that $\gm$ is a $\Pi^m_n$ reflecting cardinal. Let
$e:(\Phi,\psi,a)\rarrow \gk$ be an enumeration
of triples $(\Phi,\psi,a)$ s.t. $\Phi$ is a $\Pi^m_n$
sentence, $\psi(x,a)$ is
a first order formula in the free variable $x$,
$a \in V_{\gk}$. Also assume without loss of generality that $a\in V_{e(\Phi,\psi,a)}$ and
$e(\Phi,\psi,a)<rank(a)^{\flat}$, where $\ga^{\flat}$ is the
least inaccessible above $\ga$, and that $e$ is $1-1$
and onto $\gk$.
For each triple $(\Phi,\psi,a)$ define $g(e(\Phi,\psi,a))$ by:
\beq
 g(e(\Phi,\psi,a)) = \left\{
\begin{array}{l}
\hbox{\ The least $\gr \in C$ such that $e(\Phi,\psi,a)<\gr$ and } \\
\left(V_{\gr},\gre,S\right)\models (x\in S
\hbox{\ iff \ } \psi(x,a)) \wedge \Phi,\\ \hbox{\ if there is such a \ } 
\gr.\\ 
\ \\
e(\Phi,\psi,a)+1 \hbox{\ otherwise}.
\end{array}
\right.
\eeq
$g:\gk \rarrow \gk$ satisfies for every $\ga$, $\ga<g(\ga)$
by its definition.
Since $\gk$ is Mahlo the set of inaccessibles is stationary so
there is an inaccessible $\gm \in C$ s.t.
$g''\gm \subseteq \gm$. We claim that $\gm$ is $\Pi^m_n$
reflecting. By definition $\gm$ is inaccessible. Suppose
that $S \subseteq V_{\gm}$ is defined by
$x \in S \hbox{\ iff \ } V_{\gk} \models \psi(x,a)$.
Assume that $\left(V_{\gm},\gre,S\right) \models \Phi(S)$.
Since $\gm \in C$ there is a $\gr \in C$ s.t.
$\left(V_{\gr},\gre,S \cap \gr\right) \models \Phi$ and $S\cap \gr$ is
defined by $\psi(x,a)$. Hence $g(e(\Phi,\psi,a))$ is defined
by the first case, since $g''\gm \subseteq \gm$ and
$e(\Phi,\psi,a)<\gm$. Now by the definition of $e$,
we obtain
that $g(e(\Phi,\psi,a))<\gm$. Hence $\Phi$ reflects to
$\left (V_{g(e(\Phi,\psi,a))},\gre,S\cap g(e(\Phi,\psi,a))\right)$.
\eprf

The following theorem is the analog of the Hanf-Scott theorem on the
indescribability of weakly
compact cardinals.
\bdf
$\gk$ has the extension property iff for every $n$ and for every
$A \subseteq V_{\gk}$ first order definable over $V_{\gk}$ with
parameters from $V_{\gk}$, there is a transitive ser $X$, and
$A^X \subseteq X$ such that $\gk \in X$
and $\left(V_{\gk},\gre,A\right) \prec_n \left(X,\gre,A^X\right)$.
\edf
{\bf Remark.} One can show that $\gk$ has the extension property if and
only if for every $n$ there is a transitive $\gS_n$-elementary end
extension of $V_{\gk}$. containing $\gk$, that is, $\gk \in X$ and
$\left(V_{\gk},\gre\right) \prec_n \left(X,\gre\right)$.

Note that the fact that $\gk$ has the extension property is $\gS^1_1$.
\bthm
\label{HanfScott}
A cardinal $\gk$ is $\Pi^1_1$ reflecting iff $\gk$ is inaccessible and
has the extension property.
\ethm
The proof is identical to the proof of the Hanf-Scott
theorem
(see \cite{Kanamori} pp. 59-60).\\
\bprf

Assume that $\gk$ is inaccessible and has the extension property.
Let $A \subseteq V_{\gk}$ be definable in $V_{\gk}$.
Let $\Phi$ be a $\Pi^1_1$
formula, say of the form $\forall Y \phi(Y)$, where $\phi(Y)$ is first
order formula in the predicates $Y$ and $A$, and assume that
$\left(V_{\gk},\gre,A \right)  \models \Phi$.
Let $n$ be large enough so that the sentence
$\exists \ga\left(\forall Y \in V_{\ga+1}\left(
(V_{\ga},\gre,A \cap V_{\ga})\models \phi(Y)
\right)\right)$ is $\gS_n$ in $L_{\gre}(A)$,
the ambient language of $\left(V_{\gk},\gre,A \right)$.
Now let  $(X,\gre,A^X)$ be given by the extension property, that is,
\beq
\label{sigma_n}
(V_{\gk},\gre,A) \prec_n (X,\gre,A^X)
\eeq
where $X$ is transitive and $\gk \in X$. Therefore $A^X \cap V_{\gk}=A$.
Also note that $V_{\gk+1}^X \subseteq V_{\gk+1}$ and
$V_{\gk}^X = V_{\gk}$ (for $n$ large enough). Now since
$\left(V_{\gk},\gre,A \right)  \models \Phi $ it follows that
\[
(X,\gre,A^X) \models
\forall Y \in V_{\gk+1}\left((V_{\gk}^X,\gre,A^X \cap V_{\gk}^X)
\models \phi(Y)\right)
\]
and thus
\[
(X,\gre,A^X) \models \exists \ga \left(
\forall Y \in V_{\ga+1}\left((V_{\ga},\gre,A^X \cap V_{\ga})
\models \phi(Y)\right)\right).
\]
Hence by (\ref{sigma_n})
\[
(V_{\gk},\gre,A) \models \exists
\ga((V_{\ga},\gre,A \cap V_{\ga})  \models \Phi).
\]
Hence
\[
(V_{\ga},\gre,A \cap V_{\ga})  \models \Phi.
\]
Thus $\gk$ is $\Pi^1_1$ reflecting.

For the other direction let $\gk$ be a $\Pi^1_1$ reflecting cardinal, and
fix some $n<\go$.
Let $\gs$ be a $\Pi^1_1$ formula expressing the failure of the extension
property relative to $n$, that is there is no transitive
$\gS_n$-elementary extension $X$ of $V_{\gk}$ containing $\gk$
(see the above remark).
Let $\tau$ be a $\Pi^1_1$ formula expressing the inaccessibility of
$\gk$. Let
\[
C=\{\ga<\gk | (V_{\ga},\gre) \prec_n (V_{\gk},\gre)\}.
\]
$C$ is a club subset of $\gk$ which is $\gS_{n+1}$ definable in $V_{\gk}$ . Let $\psi$ be
the formula expressing the fact that $C$ is a club.
Therefore,
\[
(V_{\gk},\gre,C) \models \gs \wedge \tau \wedge \psi.
\]
By $\Pi^1_1$ reflection there exists an $\ga$ such that
\[
(V_{\ga},\gre,C \cap \ga) \models \gs \wedge \tau \wedge \psi.
\]
Hence $\ga$ is an inaccessible cardinal below $\gk$ which is a limit
point of $C$.
Hence $\ga \in C$, so $(V_{\ga},\gre) \prec_n (V_{\gk},\gre)$.
Let $X_0=V_{\ga} \cup \{V_{\ga} \}$ and construct an elementary
submodel $(X',E) \prec (V_{\gk},\gre)$, containing $X_0$ of cardinality
$\ga$. Let $(X',E)$ be the Skolem hull of $X_0$ inside $(V_{\gk},\gre)$.
Let $(X,\gre)$ be the transitive collapse of $(X',E)$.
It follows that $(V_{\ga},\gre) \prec_n (X,\gre)$, and $V_{\ga}\in X$
hence
\[
(V_{\ga},\gre) \models \neg \gs
\]
and this is a contradiction.
\eprf

We finish this section with the following observation on
$\Pi^1_1$ reflecting cardinals:
\begin{lemma}
\label{branch_lemma}
Let $\gk$ be a $\Pi^1_1$ reflecting cardinal. Let $T$ be a
$\gk-$tree definable in the very wide sense, then $T$ has a cofinal
branch.
\label{br_exist}
\end{lemma}
\bprf
Note that since $\gk$ is inaccessible $V_{\gk}=H_{\gk}$.
Let $T$ be a definable $\gk$-tree definable in the very wide sense.
Let $n$ be large enough
such that the assertion ``$\forall \ga \ T \hbox{\ has a cofinal branch
of length  $\ga$.}$" is
$\bgS_n$ over $\left(V_{\gk},\gre, T\right)$.
By theorem \ref{HanfScott} there is a transitive structure
$\left(X,\gre,T\right)$ such that
\beq
(V_{\gk},\gre,T) \prec_n (X,\gre,T^X).
\eeq
Since (for $n$ large enough) $V_{\gk}^X=V_{\gk} \in X$, it follows that
$T^X \cap V_{\gk}^X=T$.
Now, since $T$ is a $\gk$-tree it follows that
\[
\left(V_{\gk},\gre,T \right) \models \forall \ga \ T \hbox{\ has a cofinal
branch of length  $\ga$.}
\]
Therefore,
\[
\left(X,\gre,T^X \right) \models \forall \ga \ T^X \hbox{\ has a
 branch of length  $\ga$.}
\]
Since $\gk \in X$ we see that
$\left(X,\gre,T^X \right) \models T^X
\hbox{\ has a branch $b$  of length \ } \gk$. Because
$T^X \cap V_{\gk}=T$, this branch $b$ is really a cofinal branch through
$T$.
\eprf
\section{The forcing construction}
\label{Forcing_model}
In this section we describe the forcing construction and prove that the 
extended model
satisfies the tree property for $\go_1-$trees first order definable over
$\Hom$.
Let $\gk$ be a $\Pi^1_1$ reflecting cardinal in $V$.
Let
\[
\mP=\coll(\go,<\gk)
\]
be the Levy collapse of $\gk$ to $\go_1$
and for every $\ga<\gk$ let
\beq
\mP_{\ga}=\coll(\go, <\ga)
\eeq
be an initial segments of the forcing.
Let $G$ be a $\mP$ generic filter. Our main theorem is
\bthm
\label{tree_theorem}
$V[G] \models$ ``every definable $\go_1-$tree T has a cofinal branch''.
Moreover if $V=L$ then $V[G] \models$ ``every $\go_1-$tree T
definable in the wide sense over $\Hom$ has a cofinal branch''.
\ethm

Define the following definable partition relation :
\bdf
$ \ha_1 \stackrel{def}{\rightarrow} (\ha_1)^m_{\ga}$ iff
every partition of $ [\ha_1]^m$ into $\ga$ sets which is first order
definable over $\Hom$ (with parameters in $H_{\go_1}$), has a
homogeneous set of size $\ha_1$.
\edf
In order to prove theorem \ref{tree_theorem} we shall first prove
a definable Ramsey theorem in $V[G]$.
\blm
$V[G] \models \ha_1 \stackrel{def}{\rightarrow} (\ha_1)^2_2$
\label{part_lemma}
\elm

\bprf
Let $F: [\ha_1]^2 \rightarrow 2$ be a definable function in $V[G]$,
defined by
\beq
F(\{\ga,\gb\})=i \iff \Phi(x,\ga,\gb,i)
\eeq
where $\Phi$ is a $\bgS_n$ formula relativized to $(H_{\go_1})^{V[G]}$,
and the parameter $x$ can be taken as a real, that is,
a function $x:\go \rarrow \go$.
By the $\gk-$c.c. of the Levy collapse for every real
$x \in V[G]$ there is an $\gre<\gk$ such that  $x \in V[G_{\gre}]$, where
$G_{\gre}$ is an initial segment of the generic, which is generic for
$\mP_{\gre}=\coll(\go, < \gre)$.
Moreover by a result of Solovay (see \cite{Kanamori} proposition 10.21)
$V[G]=V[G_{\gre}][H]$, and $H$ is generic for the Levy collapse
$\coll(\go,<\gk)$. By the homogeneity of the Levy collapse every set of
ordinals definable in $V[G]$ with parameters in $V[G_{\gre}]$
is definable in $V[G_{\gre}]$, and for every formula $\Psi$ we can compute
another formula $\Phi$ such that:
\beq
V[G] \models \Psi(x,\ga,\gb,i) \iff V[G_{\gre}] \models
\Phi(x,\ga,\gb,i)
\eeq
for all $\ga,\gb<\gk$ and $i \in \{0,1\}$.
Fix such an $\gre$. Let $\tlgs$ be a $\mPeps$ name for $x$.
Now we define a $\gk-$tree $T$ which is definable
in $V_{\gk}$. For each $\ga<\gk$ let
\beq
\begin{array}{cc}
\tlh \in T_{\ga} \iff & \\
 & \tlh \hbox{\ is a \ } \mPeps\hbox{\ name for a function
from $\ga$ to \ } \{0,1\} \hbox{\ and \ } \\
 & \exists \gm_0 \forall \gm
\exists \gd>\gm  \forall \gb<\ga \forall p \in \mPeps \\
 & p \force_{\mPeps} \Phi(\tlgs,\gb,\gm_0,\tlh(\gb)) \iff  p \force_{\mPeps}
\Phi(\tlgs,\gb,\gd,\tlh(\gb))
\end{array}
\label{Tdef}
\eeq
and define $T=\cup_{\ga<\gk}T_{\ga}$. We shall write $\tlh_{\ga}$ to
denote that
$\tlh_{\ga} \in T_{\ga}$. \\
The ordering of $T$ is
\beq
\tlh_{\ga} \le \tlh_{\gb} \iff  \force_{\mPeps} \tlh_{\ga} \subset \tlh_{\gb}
\eeq
The tree $T$ is a $\gk-$tree. To prove this first
observe that for every $\ga<\gk$ $T_{\ga} \ne \emptyset$
since for every $\ga$ there are at most
$2^{2^{\ga+|\hbox{\mPeps}|}}$ many $\mPeps$ names for such functions,
and therefore
we are partitioning $\gk$ into less than $\gk$ many subsets according  to
the possible values of $\left<F(\{\gb,\gd\}) : \gb < \ga\right>$.
Secondly $|T_{\ga}| \le 2^{\ga} < \gk$.
$T$ is definable in the very wide sense in $H_{\gk}$ by
(\ref{Tdef}).
Therefore by lemma (\ref{br_exist}) $T$ has a cofinal
branch $\left<\tlh_{\ga} : \ga <\go_1 \right>$.
Work now in $V[G]$.
Let $\tlh_{\ga}(G_{\gre})$ denote the realization of the name
$\tlh_{\ga}$ in $V[G_{\gre}]$.
Let $h=\cup_{\ga<\gk} \tlh_{\ga}(G_{\gre})$, h is a function from $\gk=\ha_1^{V[G]}$ to
$\{0,1\}$. Define
\beq
A_{\ga}=\{\ga<\grg<\ha_1 | \forall \gb<\ga F(\{\gb,\grg\})=h(\gb)\}
\eeq
then for every $\ga$, $|A_{\ga}|=\ha_1$ by the definition of $\tlh_{\ga}$.
Moreover $\left<A_{\ga} : \ga<\ha_1\right>$ is a decreasing sequence of
sets. We construct $H_0$ by induction on $\ga<\ha_1$. Let $\gb_0=0$.
For each $\ga$ let $\grg_{\ga}=sup \left \{ \gb_i | i<\ga \right \}$.
Let $\gb_{\ga}=\min A_{\grg_{\ga}}$.
Let $H_0=\left\{\gb_{\ga} | \ga<\ha_1 \right\}$. By the definition of $H_0$, for
every $\ga<\gb \in H_0$ we have $F(\{\ga,\gb\})=h(\ga)$. Let
$l$ be minimal such that $|h^{-1}(l) \cap H_0|=\ha_1$. Now
$H=\left\{\ga \in H_0 | h(\ga)=l \right \}$ is the homogeneous set.
\eprf

Note that the proof really gives the following consequence \\
\[
V[G] \models \ha_1 \stackrel{def}{\rightarrow} (\ha_1)^2_{\ha_0}.
\]
If $V=L$, or there is a definable well-ordering of
$H_{\gk}$ in $V$, then the same proof yields the following:
\blm
\label{wide_ramsey}
Assume $V=L$. Let $h,A \in V[G]$ be such that $|A|=\ha_1$, and $h:[A]^2 \rarrow 2$
is a partition of $[A]^2$ into two
parts, where
both $A$ and $h$ are first order definable
(with parameters) over $\Hom^{V[G]}$.
Then there is a $B \subseteq A$  homogeneous for $h$
and $|B|=\ha_1$.
\elm
To derive theorem \ref{tree_theorem} we follow the proof that a weakly
compact cardinal
has the tree property (see lemma 29.6 of \cite{Jech}), replacing the partition property, of
a weakly compact cardinal with the definable partition lemma
(\ref{part_lemma}).\\
{\bf proof of theorem \ref{tree_theorem}}\\
Let $T=(\ha_1,<_T)$ be a definable tree on $\ha_1$, i.e. there is a
formula $\Psi(\ga,\gb,z)$ such that
\beq
\ga<_T \gb \iff  \Hom \models \Psi(\ga,\gb,z).
\eeq
We extend the partial tree ordering $<_T$ into a linear ordering as
follows:

$ \underline{\ga\prec \gb}$ iff
\begin{itemize}
\item[(i)]$\ga<_T\gb$
or
\item[(ii)]$\ga,\gb$ are $<_T$ incomparable  and if $\gz$ is the first level where the
predecessors of $\ga,\gb$, $\ga_{\gz},\gb_{\gz}$ are distinct then
$\ga_{\gz}<\gb_{\gz}$.
\end{itemize}
$\prec$ is first order definable in $\Hom$ using the definition of
$<_T$.
Now define a partition of $[\ha_1]^2$ by
\beq
F(\{\ga,\gb\})=1 \iff \ga< \gb \hbox{\ agrees with \ } \ga \prec \gb.
\eeq
Since both $\prec$ and $<_T$ are definable $F$ is definable as well.
Hence there is $H \subset \ha_1$ which is homogenous for $F$, with
$|H|=\ha_1$. Let
\beq
B=\left\{x \in \ha_1 | \ |\{\ga\in H| x<_T\ga\}|=\ha_1\right\}.
\eeq
Since every level is countable there are members of $B$ of every level.
If we'll prove that any two members of $B$ are $<_T$ comparable, then
$B$ will
be the $\ha_1$-branch. Let $x,y \in B$ be $<_T$ incomparable elements.
Assume, without loss of generality, that
$x \prec y$. Since both $x,y$ have $\ha_1$ many $<_T$ succesors in $H$
we can
find $\ga,\gb,\gm \in H$ such that $\ga<_T\gb<_T\gm$, $x<_T \ga,\gm$ and
$y<_T\gb$.
By the definition of $\prec$ we get $\ga \prec \gb$ and $\gm \prec \gb$.
Thus $F(\{\ga,\gb\})=1$ and $F(\{\gb,\gm\})=0$ , contradicting the
fact that $H$ is homogeneous for $F$.
Finally note that if we force over $L$ theorem \ref{tree_theorem} can
be strengthened to trees definable in the very wide sense. The proof is
identical using lemma \ref{wide_ramsey} instead of \ref{part_lemma}.
\eprf
\section{The lower bound}
\label{lower_bound}
In this section we prove that the definable tree property implies the
consistency of a $\Pi^1_1$ reflecting cardinal. In this section
let $\ha_1$ denote $\ha_1^V$.
\bthm
If $\ha_1$ has the definable tree property
then
\[
L \models \ha_1 \hbox{\ is a $\Pi^1_1$ reflecting cardinal}.
\]
\ethm
First we prove that
\beq
\label{inacc}
L\models \ha_1 \hbox{\ is inaccessible}.
\eeq
Assume that $\ha_1$ is not inaccessible in $L$ then there is an
$x \in \go^{\go}$ such that $\ha_1=\ha_1^{L[x]}$ (see \cite{Kanamori}
proposition 11.5). However inside
$L[x]$ there is a special Aronszajn tree $T$ which is definable from the
well ordering
of $\Hom^{L[x]}$, which is itself $\bgS_1$ definable over $\Hom^{L[x]}$
\cite{Dev}.
However $T$ cannot have a cofinal branch in $V$ since this implies
$\ha_1^{L[x]}<\ha_1$. Similarly, by relativization, for every real $x$
$\ha_1$ is inaccessible in $L[x]$.

Next we prove that $\ha_1$ is a $\Pi^1_1-$reflecting cardinal in $L$.
The proof is based on an idea from \cite{MgShSt}. We define a tree using
the $\bgS_n$ definable power set of $\ha_1$.
Note that since $\ha_1$ is inaccessible in $L$, by (\ref{inacc}), we have
$(H_{\ha_1^V})^L=L_{\ha_1^V}=\left(V_{\ha_1^V}\right)^L$.
From a cofinal branch in the tree we define an
ultrafilter on the ${\bgS_n}$ definable subsets of $\ha_1$, and
construct an ``ultrapower'' of $L_{\ha_1}$, using only functions
$\bgS_n$ definable over $\Hom$.
Note that this is a definable tree in the strong sense since
there is a definable well-oredering of the underlying set.
Let $\left< A_{\ga} | \ga< \ha_1 \right>$ be a definable enumeration of
$P^L(\ha_1) \cap \bgS_n$ of order type $\ha_1$.
Define a tree $T$ of functions by
\[
f \in T \iff f:\tau \rightarrow \{ 0,1\}, \tau<\ha_1 \hbox{\ and \ }
|\cap_{\ga<\tau} A_{\ga}^{f(\ga)}|=\ha_1
\]
where
$A^0=A$, and $A^1=\ha_1 \backslash A$.
The ordering on $T$ is $f <_T g$ iff $f \subseteq g$.
Since $\ha_1$ is inaccessible in $L$ the tree is an $\go_1-$tree.
Since the truth of $\gS_n$ formulas is a $\gS_{n+1}$ definable,
the tree $T$
is $\bgS_k$ definable over $\Hom$, for some $k$. Hence by the definable
tree property it has a cofinal branch producing a function
$b:\ha_1 \rightarrow \{0,1 \}$. $b$
defines an ultrafilter $U$ on $P^L(\ha_1) \cap \bgS_n$ by
\[
A_{\ga} \in U \iff b(\ga)=0
\]
and this ultrafilter is countably complete on $P^L(\ha_1) \cap \bgS_n$.
The ``ultrapower'' is now defined by
\beq
\begin{array}{ll}
f \in ult(L_{\ha_1},U) \iff & f:\ha_1 \rightarrow L_{\ha_1}, f \in L \\
                            &\hbox{ \ and $f$ is $\bgS_n$ definable over
$\Hom$}
\end{array}
\eeq
\beq
f \equiv g \iff \{\ga | f(\ga)=g(\ga)\} \in U
\eeq
and
\beq
f E g \iff \{\ga | f(\ga) \in g(\ga)\} \in U.
\eeq
The ``ultrapower'' is wellfounded by the completeness of the ultrafilter
and there is a $\gS_n$ embedding
$j:L_{\ha_1}\prec_n ult(L_{\ha_1},U)$, by the proof of $\L$os theorem.
However since $L_{\ha_1} \models V=L$
if $n$ is large enough $ult(L_{\ha_1},U) \models V=L$, and hence its
transitive collapse is really $L_{\ga}$ for some $\ga>\ha_1$.
Therfore $L_{\ha_1} \prec_n L_{\ga}$.
This gives the desired extension
property, and by the equivalence of the extension property and $\Pi^1_1$
reflection $\ha_1$ is $\Pi^1_1$ reflecting in $L$.
Finally by relativizing we obtain that for every real $x$,
$\ha_1$ is $\Pi^1_1$ reflecting in $L[x]$. \eprf

\section{The definable tree property and Martin's axiom}
\label{DTPMA}
In this section we investigate the consistency strengh of the Definable
Tree Property on $\ha_1$ together with Martin's Axiom and large continuum.
The exposition is based on the Shelah-Harrington paper \cite{HaSh}.
The main theorem is the following :
\bthm
Let $\ha_0<\gl$ be a cardinal satisfying $\gl^{<\gl}=\gl$. The following
are equiconsistent :
\bds
\item[1.] The definable tree property on $\ha_1 \ + \  MA \ +
 \ 2^{\ha_0}=\gl$
\item[2.] $\ha_1$ is weakly compact cardinal in $L$.
\eds
\ethm
For the $2 \Rightarrow 1$ direction, we use lemma \ref{inacc}
which proves that
the definable tree property implies that
for every real $x$ $\ha_1^{L[x]}<\ha_1$.
Now we finish by the following result from \cite{HaSh}:
\bthm
Assume MA then either there exists a real $x$ such that
$\ha_1^{L[x]}=\ha_1$, or $\ha_1$ is weakly compact in $L$.
\ethm
The proof of the other direction follows closely Kunen's forcing for
the consistency of ``$MA+$ Every set in L(R) is Lebesgue
measureable $+2^{\ha_0}=\gl$".
The model we use is Kunen's model, and we only prove that the definable
tree property holds in that model.
 For completeness we present the full proofs
of Kunen's basic observations, as presented in \cite{HaSh}.
Let $\gk$ be a weakly compact cardinal.
\blm
\label{kunen_a}
If $B$ is a complete Boolean Algebra with the $\gk-$c.c., and if
$X \subseteq B$ satisfies $|X|<\gk$, then there is a complete subalgebra
$\bar B$, s.t. $X \subseteq \bar B$ and $|\bar B|<\gk$.
\elm
\bprf
Since $B$ satisfies the $\gk-$c.c., and $\gk^{<\gk}=\gk$, there is $B'$ a
complete subalgebra of $B$ such that $X \subseteq B'$ and $|B'|\le\gk$.
Without loss of generality assume that $B' \subseteq \gk$. Let $D$ be the set
of maximal antichains of $B'$. $D \subseteq [\gk]^{<\gk}$. By $\Pi^1_1$
reflection there is an $\ga<\gk$ s.t. $B' \cap \ga$ is $<\ga$ complete
and $D \cap [\ga]^{<\ga}$ is the set of maximal antichains of
$B' \cap \ga$. Therefore $B' \cap \ga$ is a complete subalgebra of $B$.
If we choose $\ga>\sup(X)$ then $X \subseteq B'\cap \ga$.
\eprf
\blm
\label{kunen_b}
If $P_0, P_1$ are two complete Boolean algebras with $\gk-$c.c., then
$P_0 \times P_1$ has the $\gk-$c.c.
\elm
\bprf
Let $\left< <p^0_{\ga}, p^1_{\ga}> | \ga<\gk \right>$ be a sequence of
elements of $P_0 \times P_1$. Define $F:[\gk]^2 \rarrow 2 \times2$.
$F(\{\ga,\gb\})(i)=0 \hbox{\ iff\ } p^i_{\ga}, p^i_{\gb}$ are compatible.
A size $\gk$ homogeneous set for $F$ gives a size $\gk$ set of pairwise
compatible elements, since by the $<\gk-$c.c of $P_i$ the homogeneous color
is $<0,0>$.
\eprf

Assume that $\gn<\gl=\gl^{<\gl}$.
To obtain the model we iterate $\gl$ many times $\gk-$c.c.
forcings of size $<\gl$ using finite support, the same way
this is done for $MA$ (\cite{SoTn}). By lemma (\ref{kunen_b})
we can assume,
without loss of generality,
that each forcing appears $\gl$ many times in the iteration.
Denote the iteration by $B$.
To prove that every $\ha_1$ tree which is ordinal definable from a real
has a branch in $V^B$, we first prove that Ramsey theorem for $\ha_1$
holds for $L(R)$ partitions. Then we finish the same way as in the proof
of theorem (\ref{tree_theorem}).

Let $F:[\ha_1]^2 \rarrow 2$ be definable in $V^B$ from a real $x$, and
ordinal parameters. By lemma (\ref{kunen_a}) $x$ is generic for a countable subalgebra $\bar B$.
Moreover since each forcing in the iteration appears unboundedly many
times  we can assume , without loss of generality, that $B/{\bar B}$ is an homogeneous
iteration over $V^{\bar B}$ which is built the same way as $B$ is
built over $V$. Therefore every value of $F$ is decided inside
$V^{\bar B}$. Since $|{\bar B}|<\gk$ we can construct (in $V$) a $\gk$
tree of $\bar B$ names
of possible values of $F$, the same way as we have built it in lemma
(\ref{part_lemma}).
Now by the weak compactness of $\gk$ we can find a branch through
that tree. Finally we use the branch to obtain the homogeneous set the
same as we have done in lemma (\ref{part_lemma}).

Note that since $\gk$ is weakly compact we are able to find homogeneous
sets for every $L(R)$ coloring, and not just first order definable over
$H_{\go_1}$.
\eprf

Finally we remark that like in Solovay's proof of Lebesgue
measurability \cite{solovay70}, just adding the negation of the Continuum
Hypothesis does not add to the consistency strength. We have
to notice that the product of collapsing $\gk$ to $\go_1$
and then adding $\gl$ Cohen reals satisfies the $\gk-$c.c. The product is 
also homogeneous enough.
Moreover every real belongs to a small generic
extension, which is complemented by a homogeneous forcing.
Hence the same situation as in theorem \ref{tree_theorem} still holds,
and the argument there can be carried out.

\section{The definable tree property on higher cardinals}
\label{DTPGCH}
This section contains several remarks regarding the
definable tree
property on cardinals above $\ha_1$. By Mitchel's and Silver's
results \cite{Mi}, the tree property on $\ha_2$ is equiconsistent
with the existence of a weakly compact cardinal.
However it is well known that assuming GCH, or even
$\gl^{<\gl}=\gl$, we have a $\gl^+$ special Aronszajn tree.
This is proved by the
same construction as Aronszajn's construction on
$\ha_1$, using the fact that under this hypothesis
there is a universal linear order of cardinality $\gl$.

We will prove that assuming a $\Pi^1_1$ reflecting cardinal
the definable tree property on $\ha_2$
is consistent with GCH.
Then we consider the property of having the definable tree
property on succesive cardinals. We will generalize our forcing argument
to prove that if it is consistent that there are $\go$ $\gP^1_1$
reflecting cardinals then it is consistent to have the definable
tree property on $\left< \ha_i : 1 \le i <\go \right>$.
\bthm
Assume $\gk$ is a $\gP^1_1$ reflecting cardinal in L. Let $P=Coll(\ha_1, <\gk)$.
Let $G$ be $P$ generic over L. Then
\[
L[G] \models \hbox{GCH and \ } \ha_2 \hbox{\ has the definable tree property.}
\]
\ethm
\bprf
The proof is identical to the $\ha_1$ case.
Using the homogeneity of the Levy collapse, and the fact
that $\gk$ is $\gP^1_1$ reflecting. We just build a tree of
possible values for the definable function. Then using the
definabiltiy of the function we obtain definability of the tree
and thus we can use the branch to prove the definable version
of Ramsey theorem for $\ha_2$. The fact that GCH holds in $V[G]$, is
proved in the usual way, by assuming $GCH$ in the ground model.
\eprf

To obtain the definable tree property on all $\ha_n$'s from
a sequence of $\go$ $\gP^1_1$ reflecting cardinals, $\gk_0<\gk_1<\ldots$
just
iterate the forcings $Coll(\ha_n,<\gk_n)$ with finite
support. The proof that
the definable tree property holds for every n, is a
straightforward generalization and will not be given here.
\section{Acknowledgement}
I would like to thank Menachem Magidor for helpful
discussions, and the anonymous referee for his comments
which greatly improved the presentation of the results.

\end{document}